# OPTIMAL PROPERTIES OF CENTROID-BASED CLASSIFIERS FOR VERY HIGH-DIMENSIONAL DATA

By Peter Hall and Tung Pham

*University of Melbourne*

We show that scale-adjusted versions of the centroid-based classifier enjoys optimal properties when used to discriminate between two very high-dimensional populations where the principal differences are in location. The scale adjustment removes the tendency of scale differences to confound differences in means. Certain other distance-based methods, for example, those founded on nearest-neighbor distance, do not have optimal performance in the sense that we propose. Our results permit varying degrees of sparsity and signal strength to be treated, and require only mild conditions on dependence of vector components. Additionally, we permit the marginal distributions of vector components to vary extensively. In addition to providing theory we explore numerical properties of a centroid-based classifier, and show that these features reflect theoretical accounts of performance.

## 1. Introduction.

1.1. *Motivation and summary.* Suppose we observe samples $\mathcal{X}$ and $\mathcal{Y}$, both consisting of $p$-vectors, drawn by sampling randomly from respective populations $\Pi_X$ and $\Pi_Y$. In this paper we establish optimality properties for classifiers based on the centroid method in cases where $p$ is large and sample sizes are, generally, much smaller. For the applications we have in mind, sample sizes can be quite small indeed; for example, in genomic problems $p$ is typically in the thousands or tens of thousands, but training sample sizes may be only in the teens, or even less. It is shown that in cases such as this, a scale-adjusted version of the classifier is able to discriminate in an optimal way between populations that differ in terms of location. Scale adjustment removes the tendency for scale to confound location differences when using distance-based classifiers, and permits the method to enjoy high levels of performance when location differences are relatively small.

In order to outline our main results, let us suppose that location differences are present in a proportion $q$ of the $p$ components; that both training sample sizes are at least as large as 2, and are of similar magnitude, $\nu$ say and that the components of the data vectors are not too strongly correlated, in particular that the maximum of the sum of absolute values of covariances, against any particular component, is bounded. Then a good classifier can correctly distinguish between the populations that correspond to the training samples, provided that the size of the location differences is a sufficiently large constant multiple of $(\nu p q^2)^{-1/4}$. Moreover, in minimax terms this size of distance is the minimum possible for accurate discrimination.

These results hold for large values of the dimension, $p$, and in particular they are valid in cases where dimension is of larger order than the training-sample sizes. However, the results can fail if sample sizes are of a larger order than $p$, for example, if $p$ is held fixed while samples increase. Therefore, our results specifically address the case where dimension is high.

In our lower-bound analysis we impose the condition that $q$ exceeds a constant multiple of $(\nu/p)^{1/2}$, thereby preventing sparsity, indexed by $q$, from being too low. This assumption implies that $(\nu p q^2)^{-1/4}$ is bounded above by a constant multiple of $\nu^{-1/2}$, and entails boundedness of the location differences. However, our work does not require $\nu$, denoting the order of magnitude of training-sample size, to diverge; $\nu$ can be held fixed, although it can be chosen to diverge if desired. Therefore our results encompass cases where the location differences are bounded away from zero as $p$ increases, as well as instances where the differences converge to zero.

1.2. *Interpretation.* First we note that results of the type discussed above hold only in the very high-dimensional cases that are the subject of our work, and not in more conventional settings. To indicate why, let us simplify matters by taking $q = 1$. In this setting it is readily shown that if $p$ is held fixed, but $\nu$ is permitted to increase, then simple distance-based classifiers can detect location differences that are of order $\nu^{-1/2}$ in size. However, fixing $p$ and varying $\nu$ in the convergence-rate formula $(\nu p q^2)^{-1/4} = (\nu p)^{-1/4}$ would suggest, incorrectly, that the best rate is only $\nu^{-1/4}$. Therefore the formula is not applicable to cases where dimension is much smaller than sample size. More specifically, the fact that the critical quantity $(\nu p q^2)^{-1/4}$ involves the exponent $-\frac{1}{4}$, rather than $-\frac{1}{2}$ which arises in more conventional settings, underscores the challenge of undertaking classification using small samples of high-dimensional data, rather than large samples of low-dimensional data.

Among classification problems that are relatively difficult to solve are those where the location differences that distinguish the two populations are so irregular as to resemble stochastic processes. In such cases, classifiers can readily confuse location differences with additive random noise. Therefore, when establishing lower bounds we interpret location differences as random variables that have the same distribution (after rescaling) as the noise. Our upper-bound results also permit this treatment.

1.3. *Comparison with other classifiers.* Other classifiers, for example, based on nonparametric function estimation or $k$-nearest neighbor methods, are competitive under suitable conditions. Several classifiers can be interpreted, either explicitly or implicitly, as empirical approximations to the Bayes classifier. For example, Stone (1977) discusses empirical classifiers based on function approximations, and Cover (1968) and Devroye and Wagner (1982) address $k$-nearest-neighbor methods. Using the latter approach, and in low-dimensional settings, if $k$ is chosen to diverge appropriately as sample size increases then the classifier can achieve the same first-order asymptotic performance as the Bayes method. This is achieved through the classifier implicitly estimating the unknown densities, $f_X$ and $f_Y$, say, of the two populations, and using them in a manner which is first-order equivalent to the Bayes rule, that is, assigning a new data value, $Z$, to $\Pi_X$ if $f_X(Z) > f_Y(Z)$ and assigning $Z$ to $\Pi_Y$ otherwise. The empirical approaches suggested by Stone (1977) and Hall and Kang (2005) do this more explicitly. If $p$ increases sufficiently slowly as the training sample sizes diverge then empirical classifiers such as these can strongly outperform the centroid-based method.

However, both explicit and implicit estimation of $f_X$ and $f_Y$ are ineffective when the dimension is of the same order as, or of larger than, the sample sizes. There, methods such as the centroid-based classifier and the support vector machine come into their own. Both these methods exhibit the optimal performance expressed by Theorems 1 and 2. In the case of the support vector machine we need somewhat more restrictive conditions than those that we impose in Section 3, and in particular which require the training sample sizes to diverge no more quickly than $p^{1/10}$. A proof in that case is given in the unpublished Ph.D. thesis of the second author.

1.4. *Related work.* The literature on statistical classification is particularly extensive, and we shall provide here only a brief pointer to relatively recent literature. Hastie, Tibshirani and Friedman (2001) give a benchmark survey of statistical learning, and Dudoit, Fridlyand and Speed (2002) provide an authoritative comparison of the performance of statistical classifiers. Dabney (2005), Dabney and Storey (2005, 2007), Tibshirani et al. (2002) and Wang and Zhu (2007) discuss the application of centroid-based classifiers to genomic data. Many other contributions are written from the viewpoint of engineering, computer science and other fields, rather than statistics, and address applications in areas ranging from image analysis [e.g., Cootes et al. (1993)] and forestry [e.g., Franco-Lopez, Ek and Bauer (2001)] to speech recognition [e.g., Bilmes and Kirchhoff (2004)] and chemometrics [e.g., Schoonover, Marx and Zhang (2003)]. They include work on the development of transformation methods for improving classifier performance [e.g., Sinden and Wilfong (1992), Simard, Lecun and Denker (1993) and Wakahara, Kimura and Tomono (2001)]. Chan and Hall (2009) provide background to scale adjustment. Nearest-neighbor methods are discussed by Dasarathy (1990) and Shakhnarovich, Darrell and Indyk (2005). Van der Walt and Barnard (2006) give a recent account of classifier performance. Duda, Hart and Stork (2001) provide a book-length treatment of classifiers in the context of pattern recognition.

## 2. Scale adjustment.

2.1. *Scale-adjusted centroid-based classifier.* A standard centroid-based classifier can be defined as follows. Let $\mathcal{X} = \{X_1, \ldots, X_m\}$ and $\mathcal{Y} = \{Y_1, \ldots, Y_n\}$ denote random samples of $p$-vectors from populations $\Pi_X$ and $\Pi_Y$, respectively, and write $\bar{X} = m^{-1} \sum_i X_i$ and $\bar{Y} = m^{-1} \sum_j Y_j$ for the respective sample means. Put

$$T(Z) = \|Z - \bar{Y}\|^2 - \|Z - \bar{X}\|^2. \tag{2.1}$$

Given a new data vector $Z$ from one of the two populations, classify $Z$ as coming from $\pi_X$ if $T(Z) > 0$, and assign $Z$ to $\Pi_Y$ if $T(Z) \leq 0$.

This classifier is used frequently to distinguish between two populations on the basis of location differences. In that setting it enjoys good performance if the training sample sizes $m$ and $n$ are reasonably large, but in other cases its effectiveness can be hampered by excessive scale differences. A simple adjustment removes this difficulty. Specifically, define

$$\hat{\tau}_X^2 = \frac{1}{2m(m-1)} \sum_{i_1=1}^{m} \sum_{i_2=1}^{m} \sum_{k=1}^{p} (X_{i_1 k} - X_{i_2 k})^2,$$

$$\hat{\tau}_Y^2 = \frac{1}{2n(n-1)} \sum_{i_1=1}^{n} \sum_{i_2=1}^{n} \sum_{k=1}^{p} (Y_{i_1 k} - Y_{i_2 k})^2,$$

denoting unbiased estimators of

$$\tau_X^2 = \sum_{k=1}^{p} E(X_{ik} - EX_{ik})^2, \qquad \tau_Y^2 = \sum_{k=1}^{p} E(Y_{ik} - EY_{ik})^2,$$

respectively. The scale-adjusted form of $T(Z)$, whether defined by (2.1) or (2.2), is

$$T_{\text{sa}}(Z) = T(Z) + m^{-1} \hat{\tau}_X^2 - n^{-1} \hat{\tau}_Y^2. \tag{2.2}$$

Scale adjustments of other distance-based classifiers are also effective, but in general the adjustments differ from that given in (2.2).

From some viewpoints the correction at (2.2) provides an adjustment of bias, rather than scale. However, if we were to refer to it as a bias adjustment then it might be interpreted as a means of diminishing the effects of differences between the locations of populations $\Pi_X$ and $\Pi_Y$. To the contrary, it removes the effects of scale in order that location differences might be made more pronounced, rather than diminished.

The quantity $T_{\text{sa}}(Z)$ is an unbiased estimator of the signed sum of squares of distances among means

$$E\{T_{\text{sa}}(Z) \mid Z\} = s(Z) \sum_{k=1}^{p} (EX_{ik} - EY_{ik})^2,$$

where $s(Z) = 1$ if $Z$ is from $\Pi_X$, and $s(Z) = -1$ if $Z$ comes from $\Pi_Y$. Therefore, unlike $T(Z)$, the expected value of which is given by

$$E\{T(Z) \mid Z\} = s(Z) \sum_{k=1}^{p} (EX_{ik} - EY_{ik})^2 + n^{-1}\tau_Y^2 - m^{-1}\tau_X^2$$

in the centroid method approach, $T_{\text{sa}}(Z)$ focuses sharply on component-wise differences among means.

If it should happen that $m^{-1}\tau_X^2 = n^{-1}\tau_Y^2$, for example, if $m = n$ and the populations have identical average scales, then scale adjustment is not necessary. In this context our results for the classifier based on $T_{\text{sa}}(Z)$, in particular result (3.4) in Section 3.2, hold also for the standard classifier based on $T(Z)$.

2.2. *Scale adjustment in other contexts.* It can be seen from the definition of a centroid-based classifier that it endeavors to focus on differences in location, rather than in scale. It shares this feature with most other distance-based classifiers, for example, the support vector machine and distance-weighted discrimination. However, for all these methods, differences in scale can confound differences in location to such an extent that the classifier can finish up assigning $Z$ to whichever population has least variation, regardless of whether $Z$ comes from $\Pi_X$ or $\Pi_Y$.

One of the worst offenders in this regard is the standard nearest-neighbor method. If the populations $\Pi_X$ and $\Pi_Y$ have component-wise average variances equal to $\sigma_X^2$ and $\sigma_Y^2$, respectively, and component-wise average squared location differences equal $\mu^2$, then the nearest-neighbor classifier gives asymptotically correct discrimination, as $p \to \infty$, if and only if

(2.3) $$\mu^2 > |\sigma_X^2 - \sigma_Y^2|.$$

If $\mu^2 < |\sigma_X^2 - \sigma_Y^2|$ then, with probability converging to 1 as $p \to \infty$, the nearest-neighbor method assigns $Z$ to whichever of $\Pi_X$ and $\Pi_Y$ has least component-wise average variance, regardless of whether $Z$ came from $\Pi_X$ or $\Pi_Y$. In contrast to (2.3), the support vector machine and centroid-based classifiers require only

$$\mu^2 > |\sigma_X^2 m^{-1} - \sigma_Y^2 n^{-1}|, \tag{2.4}$$

where $m$ and $n$ denote the training-sample sizes for $\Pi_X$ and $\Pi_Y$, respectively. [These results hold in cases where $p$ is very large relative to $m$ and $n$, and under conditions discussed by Hall, Marron and Neeman (2005).] From (2.4) we see that, for support vector machine and centroid-based classifiers, the effects of increasing training-sample size can quickly reduce the impact of scale differences. However, in view of (2.3) this opportunity does not arise in the case of standard nearest-neighbor methods. In some problems the sample size issue is becoming less serious over time, as more data accumulate. However in other settings, for example, in the new uses of microarrays, the issue of small sample size can still be very important.

Of course, if we felt that that (2.3) or (2.4) correctly captured the ways in which location and scale worked together to jointly characterise populations $\Pi_X$ and $\Pi_Y$, then we would not introduce the scale adjustment suggested in Section 2.1. However, in practice one often feels that the differences between populations that are of interest are primarily those of location, not scale. For example, this tends to be the case with genomic data.

The measures of performance discussed above address relatively subtle properties, where the "signal" that gives rise to location differences is at least bounded, if not small. By way of contrast, some related work on classifier performance [see, e.g., Hall, Pittelkow and Ghosh (2007)] addresses instances where the signal, when it is present, is unboundedly large, and in fact diverges to infinity as $p$ increases. In such cases a scale adjustment is not necessary since the effect of uncorrected scale is of smaller order than the impact of the signal.

An alternative approach to scale adjustment is to empirically correct each component for scale before incorporating it in the classifier, in the manner of a $t$-statistic. If the scales of different components are genuinely different, for example, with some referring to weight and the others to distance, then standardisation is essential. Fortunately, in many of the applications to which classifiers are put the components have identical scales. For instance, in applications to genomic data the $j$th component of a data vector $X_i$ or $Y_i$ typically represents the extent to which the $j$th gene is differentially expressed, or "switched on," and is on the same scale for each gene.

In problems where scale standardizations is necessary, for example, to accommodate heteroscedasticity among vector components, small sample sizes can lead to problems when dividing by standard deviation estimators. These difficulties can be alleviated by using a ridge parameter or a related approach to regularisation, for example, the band-matrix inversion method of Bickel and Levina (2008).

## 3. Upper bound to classifier performance.

### 3.1. *Model for data.* We use the following data model:

(3.1) $X_{ik} = \delta a_k I_k + M_{ik}$, $Y_{jk} = \delta b_k J_k + N_{jk}$ and $Z_k = \delta c_k K_k + Q_k$ where (a) $\vec{M}_i = (M_{i1}, M_{i2}, \ldots)$, $\vec{N}_j = (N_{j1}, N_{j2}, \ldots)$ and $\vec{Q} = (Q_1, Q_2, \ldots)$ are infinite sequences of random variables with finite, zero means, (b) $\vec{M}_1, \vec{M}_2, \ldots$ are independent and identically distributed, $\vec{N}_1, \vec{N}_2, \ldots$ are independent and identically distributed and the $\vec{M}_i$'s, the $\vec{N}_j$'s and $\vec{Q}$ are independent, (c) $a_1, a_2, \ldots$ and $b_1, b_2, \ldots$ are sequences of constants and $I_1, I_2, \ldots$ and $J_1, J_2, \ldots$ are sequences of zeros and ones, (d) $\delta > 0$ is a deterministic function of $m$, $n$ and $p$, (e) $\min(m,n) \geq 2$ and (f) either $(c_k, K_k) = (a_k, I_k)$ for all $k$, or $(c_k, K_k) = (b_k, J_k)$ for all $k$.

In particular, we make no assumptions about the relationships among the noise distributions for the $X$ and $Y$ populations. For example, we do not ask that the distributions of $\vec{M}_1$, $\vec{N}_1$ and $\vec{Q}$ be related in any sense. Condition (e) is needed so that we can estimate the scale of the data; variability generally cannot be accessed empirically if either $m$ or $n$ equals 1. However, (e) is unnecessary if $m^{-1}\tau_X^2 = n^{-1}\tau_Y^2$ and we use the classifier based on $T(Z)$, rather than on $T_{\text{sa}}(Z)$. Condition (f) asserts that the pattern of the component means, $\delta c_k K_k$, for the new datum $Z$ is identical to that for either the $X$ or the $Y$ data. In particular, we describe differences between the two populations only in terms of location differences.

It might be thought that in the latter respect, the nonadjusted classifier based on $T(Z)$ enjoys potential advantages since it is influenced by differences in scale as well as differences in location. However, the nonadjusted classifier can actually be seriously misled by scale differences. See, for example, Chan and Hall (2009).

### 3.2. *Main results.* Define $\nu = \min(m, n)$. We assume that, for all $k \geq 1$, fourth moments of $M_{1k}$ and $N_{2k}$ exist, and second moments of $Q_k$ exist; and, more specifically, that the constants

(3.2) $$D_1 = \sup_{p \geq 1} \max \left[ \sup_{k_1 \geq 1} \sum_{k_2=1}^{\infty} |\text{cov}(M_{1k_1}, M_{1k_2})|, \right.$$

$$\sup_{k_1 \geq 1} \sum_{k_2=1}^{\infty} |\mathrm{cov}(N_{1k_1}, N_{1k_2})|, \sup_{k_1 \geq 1} \sum_{k_2=1}^{\infty} |\mathrm{cov}(Q_{k_1}, Q_{k_2})| \Bigg],$$

$$(3.3)\ D_2 = \sup_{p \geq 1} \max \left\{ \sup_{k_1 \geq 1} \left| \sum_{k_2=1}^{\infty} \mathrm{cov}(M_{1k_1}^2, M_{1k_2}^2) \right|, \sup_{k_1 \geq 1} \left| \sum_{k_2=1}^{\infty} \mathrm{cov}(N_{1k_1}^2, N_{1k_2}^2) \right| \right\}$$

are finite. Empirical evidence indicates that correlations among gene expression levels are often quite low, for example, in the range 0.08 to 0.01 at distances of between two and 10 base pairs, respectively [Mansilla et al. (2004), Messer and Arndt (2006)]. More generally, decay can occur at either an exponential or a reasonably fast polynomial rate [Almirantis and Provata (1999)].

This condition amounts to an assumption about the strength of dependence among the components of data vectors. To illustrate the implications of the condition we note that if the processes $\{M_{11}, \ldots, M_{1p}\}$, $\{N_{11}, \ldots, N_{1p}\}$ and $\{Q_1, \ldots, Q_p\}$ are all stationary and Gaussian, all with zero means and the same autocovariance function $\gamma(j) = \mathrm{cov}(Q_k, Q_{k+j})$, then finiteness of $D_1$ and $D_2$ is equivalent to convergence of the series $\sum_j |\gamma(j)|$. This is a mild assumption; the covariance can decay as slowly as $j^{-1-\eta}$, for any $\eta > 0$, and Theorem 1 will hold.

Define $d_k = a_k I_k - b_k J_k$, $d = (d_1, \ldots, d_p)$ and $\|d\|^2 = \sum_k d_k^2$. Let $T(Z)$ and $T_{\mathrm{sa}}(Z)$ be as at (2.1) and (2.2). In particular, $T(Z)$ is the centroid-method classifier. A proof of the following theorem is given in a longer version of this paper [Hall and Pham (2009)].

THEOREM 1. *Assume the model at (3.1), and in particular suppose that* (a)–(f) *there hold. Then there exists a constant $B > 0$, depending only on $D_1$ and $D_2$ at (3.2) and (3.3), such that*

$$(3.4) \qquad E\{T_{\mathrm{sa}}(Z) - \delta^2 s(Z) \|d\|^2\}^2 \leq B(\nu^{-1} p + \delta^2 \|d\|^2).$$

*Under the same assumptions, except that condition* (e) $\min(m, n) \geq 2$ *can now be dropped, we have instead of (3.4),*

$$(3.5)\quad E\{T(Z) - \delta^2 s(Z) \|d\|^2 - \tfrac{1}{2}(m^{-1}\tau_X^2 - n^{-1}\tau_Y^2)\}^2 \leq B(\nu^{-1} p + \delta^2 \|d\|^2).$$

3.3. *Implications for probability of correct classification.* Assume for simplicity that $I_k = J_k$ for each $k$. (The latter condition implies that the "signal" is present at the same locations in the $X$ and $Y$ populations.) Suppose too that

$$(3.6) \qquad W_1 pq \leq \|d\|^2 \leq W_2 pq, \qquad m + n \leq W_2 \min(m, n) = W_2 \nu,$$

where $0 < W_1 < W_2 < \infty$ are constants, and $q \in (0, 1]$ is an "index of sparsity." For example, if $I_k \neq 0$ for just $pq$ values of $k$, and if the sum of $(a_k - b_k)^2/pq$ over these indices is bounded away from zero and infinity, then the first part of (3.6) holds and $q$ denotes the proportion of components, in either the $X$ or $Y$ populations where the signals have an opportunity to be nonzero. Of course, we permit $q$ to vary with $p$ as the latter increases.

We also assume that $\nu \leq Cp$ where $C > 0$ is a positive constant. Therefore, the number of dimensions is at least as large as a constant multiple of sample size. Without this condition, the results that we shall describe below are generally false. For example, they fail if $p$ is held fixed as $\nu$ varies. In that setting it is readily shown that a classifier can detect alternatives distant $\nu^{-1/2}$, rather than $\nu^{-1/4}$, apart; the latter result would follow from the results given below if we were to take $p$ and $q$ fixed and permit $\nu$ to increase. These differences point to the intrinsic difficulty of undertaking classification using high-dimensional data in small samples, as distinct from low-dimensional data in large samples.

Take $\delta = c(\nu p q^2)^{-1/4}$ where $c > 0$ denotes a fixed constant. Let $\mathcal{M} = \mathcal{M}(C, D_1, D_2, W_1, W_2)$ denote the set of all models prescribed by the constraint $\nu \leq Cp$ [where it is assumed that (3.2), (3.3) and (3.6) hold for the constants $D_1, D_2, W_1, W_2$] and by conditions (a)–(f) in (3.1) [where we take $\delta = c(\nu p q^2)^{-1/4}$, for $c > 0$ fixed]. Then the following result holds [see Hall and Pham (2009) for a proof]:

COROLLARY 1. *If (3.4) and (3.6) hold then*

$$\lim_{c \to \infty} \limsup_{p \to \infty} \sup_{\text{model} \in \mathcal{M}} \{P(\text{the classifier } T_{\text{sa}} \text{ assigns } Z \text{ to } \Pi_X | Z \in \Pi_Y)$$
(3.7)
$$+ P(\text{the classifier } T_{\text{sa}} \text{ assigns } Z \text{ to } \Pi_Y | Z \in \Pi_X)\} = 0.$$

*That is, if the signals are distributed with sparsity $q$ and are of size approximately $c(\nu p q^2)^{-1/4}$, then the probability that the classifier based on $T_{\text{sa}}$ makes the incorrect decision can be rendered arbitrarily close to 0 for all sufficiently large $p$ and uniformly over all models in the class $\mathcal{M}$ by taking $c$ sufficiently large.*

Results such as (3.4), (3.5) and (3.7) all have analogues in settings where the "constants" $a_k$ and $b_k$ are interpreted as random variables. See, for example, (4.5) in Section 4.

Generally speaking, (3.7) fails if the scale adjustment suggested in Section 2.1 is not incorporated, unless $\nu$ is at least as large as a constant multiple of $p$. Indeed, it can be shown that if $|m^{-1}\tau_X^2 - n^{-1}\tau_Y^2|$ is larger than a sufficiently large constant multiple of $\delta^2 \|d\|^2$ (and this condition is often satisfied if $\nu < \text{const.} p$), then the probability of misclassification can be bounded away from zero as $p$ diverges. These results point to the desirability of including the scale adjustment when defining the classifier.

## 4. Lower bound to classifier performance.

4.1. *Data model for lower bound.* Assume we observe

$$(4.1) \quad X_{ik} = \delta A_k I_k + M_{ik}, \qquad Y_{jk} = \delta B_k J_k + N_{jk}, \qquad Z_k = \delta C_k K_k + Q_k,$$

where (i) $1 \leq i \leq m$ and $1 \leq j \leq n$; (ii) the random variables $A_k$, $B_k$, $M_{ik}$, $N_{jk}$ and $Q_k$ are normal $N(0,1)$; (iii) these variables, and $I_k$ and $J_k$, are totally independent for $1 \leq i \leq m$, $1 \leq j \leq n$ and $1 \leq k \leq p$; (iv) $I_k$ and $J_k$ are identically distributed, with $P(I_k = 0) = 1 - q$ and $P(I_k = 1) = q$, (v) $\delta > 0$ and $0 < q \leq 1$ and (vi) either $(C_k, K_k) = (A_k, I_k)$ for all $k$, or $(C_k, K_k) = (B_k, J_k)$ for all $k$. It is desired to distinguish between the two cases in (vi) using only the data at (4.1). For example, determining that $(C_k, K_k) \equiv (A_k, I_k)$ corresponds to classifying $Z_k$ as coming from the $X$ population. We permit $m$, $n$ and $q$ to depend on $p$, which we take to diverge to infinity.

By permitting $q$ to converge to zero as $p$ diverges we can ensure a degree of sparsity in the signals. However, we do not insist that $q$ becomes small as $p$ increases; for example, our assumptions permit $q$ to be held fixed, at 1, for all $p$.

Provided the likelihood-ratio statistic is asymptotically normally distributed, that quantity provides asymptotically optimal discrimination between the cases $(C_k, K_k) \equiv (A_k, I_k)$ and $(C_k, K_k) \equiv (B_k, J_k)$ in (4.1). A necessary condition for asymptotic normality is

$$(4.2) \qquad \max(m+1, n+1)\delta^2 \leq C,$$

where $C > 0$ is arbitrary but fixed. We shall make this assumption.

To indicate the implications of (4.2) we note that when this condition holds, the bias and error-about-the-mean contributions to the likelihood-ratio statistic are of sizes $\omega \equiv mpq^2\delta^4$, and $\omega^{1/2}$, respectively. Therefore, if $\omega$ is small then the bias, which reveals the difference between the cases $(C_k, K_k) \equiv (A_k, I_k)$ and $(C_k, K_k) \equiv (B_k, J_k)$, is submerged in noise, and it is impossible, even when using the likelihood-ratio method, to distinguish effectively between the cases. On the other hand, if $\omega$ is large, then the cases can be distinguished with high probability. It is in the intermediate setting, where $\omega$ is not far from 1, that classification is marginal; see Theorem 2, below. In such instances, if it should be the case that $m/(pq^2)$ diverges along a subsequence, and if $\delta = c(mpq^2)^{-1/4}$ as in Theorem 1, then $m\delta^2$ must also diverge along that subsequence, contradicting (4.2). Therefore the context of our work implies that $m/(pq^2)$ is bounded, which in turn entails a lower bound to sparsity; for a constant $C > 0$,

$$(4.3) \qquad C(m/p)^{1/2} \leq q \leq 1.$$

4.2. *Optimal convergence rates for the model at (4.1).* Write $P_X$ and $P_Y$ for probability measure under (4.1) in the respective cases $C_k \equiv A_k$ and $C_k \equiv B_k$. Let $\widehat{\chi}$ denote a measurable function of the data $X_{ik}$ (for $1 \leq i \leq m$), $Y_{jk}$ (for $1 \leq j \leq n$) and $Z_k$, all for $1 \leq k \leq p$. Let $\widehat{\chi}$, a random quantity, be a measurable function of the data $X_{ik}$, $Y_{jk}$ and $Z_k$, for $1 \leq i \leq m$, $1 \leq j \leq n$ and $1 \leq k \leq p$, and taking only the values $X$ and $Y$. In particular, $\widehat{\chi}$ can be interpreted as a classifier that ascribes $Z$ to either $\Pi_X$ or $\Pi_Y$. Write $\mathcal{C}$ for the set of all such classifiers.

The theorem below asserts that, unless $\delta$ is a relatively large constant multiple of $(mpq^2)^{-1/4}$, no classifier can effectively distinguish between the cases $(C_k, K_k) \equiv (A_k, I_k)$ and $(C_k, K_k) \equiv (B_k, J_k)$. Together with Theorem 1 it shows that the scale-adjusted classifier introduced in Section 2.1 has an asymptotically optimal ability to distinguish between the two populations.

Take $\delta$ in (4.1) to be given by $\delta = c(mpq^2)^{-1/4}$ where $c > 0$ is fixed.

THEOREM 2. *Assume the model in Section 4.1, and in particular suppose that* (i)–(vi) *there hold. Suppose too that, as p diverges, the positive integers m and n, and $q \in (0, 1]$, are such that (4.3) holds for a constant $C > 0$, and the ratio $m/n$ is bounded away from zero and infinity. Then, for all sufficiently small $c > 0$,*

$$(4.4) \qquad \liminf_{n \to \infty} \inf_{\widehat{\chi} \in \mathcal{C}} \{P_X(\widehat{\chi} = B) + P_Y(\widehat{\chi} = A)\} > 0.$$

The assumption that the "signals," represented by the terms $\delta A_k$ and $\delta B_k$ in (4.1), are random, gives them an irregular character and makes classification relatively challenging. If we take $A_k$ and $B_k$ to be fixed constants, not depending on $k$, then the classification problem is significantly simpler, and successful classification is possible for values of $\delta$ that are an order of magnitude smaller than those discussed in Theorem 2. In the model introduced at (3.1) we effectively conditioned on $A_k$ and $B_k$, treating them as constants $a_k$ and $b_k$. This is a minor alteration, however. In particular, (3.4) continues to hold if we give $a_k$ and $b_k$ the distributions of random variables, for example, as in point (ii) immediately below (4.1), and if we take expectations on both sides of (3.4). Arguing in this way the following analogue of (3.7) can be derived under the assumptions of Theorem 2.

THEOREM 3. *Assume the conditions of Theorem 2, and in particular that $\delta$ in (4.1) is defined by $\delta = c(mpq^2)^{-1/4}$. Then*

$$(4.5) \qquad \lim_{c \to \infty} \liminf_{p \to \infty} \min[P_X\{T_{\mathrm{sa}}(Z) > 0\}, P_Y\{T_{\mathrm{sa}}(Z) < 0\}] = 1.$$

Together, (4.4) and (4.5) establish optimality of the centroid-based classifier.

**5. Numerical properties.** An extensive simulation study is summarised by Hall and Pham (2009). It treats both moving-average and GARCH models Fan and Yao (2003) for the data vectors $\vec{M}$, $\vec{N}$ and $\vec{Q}$, and provides numerical evidence of theoretical properties reported in Sections 3 and 4. For example, it shows, as argued in theoretical terms in Corollary 1, that if the value of $\delta$, in the model at (3.1), is chosen so that a given, fixed percentage of classifications is correct, then $\delta$ changes with $m$ in proportion to $m^{1/4}$ if $p$ (the dimension) and $q$ (the level of sparsity) are kept fixed.

Below we report the results of sampling experiments performed using the KDD 2008 dataset. The data are available at http://www.kddcup2008.com and contain information derived from X-ray images of breast cancer patients. Two supplementary files are also provided, Features.txt and Info.txt. The Features file contains information about 102,294 suspicious regions, each described by $p = 117$ features. The Info file provides additional information about each region in the Features file. The latter file gives 11 columns describing 11 characteristics of each region. For example, the first column contains labels that indicate whether the corresponding region was malignant or benign. To simplify the classification problem we used only information about this label (i.e., malignant or benign) of each region, and ignored other information in the Info file; we used the label information only to create the samples and to assess classifier performance. Our dataset therefore contained 623 data vectors corresponding to malignant regions, and 101,671 vectors from benign regions $(623 + 101{,}671 = 102{,}294)$.

We used the KDD data to compare five methods: scale-adjusted versions of the nearest neighbor; (NN) support vector machine (SVM) and centroid-based classifiers; the scaled variance (SV) classifier for which the analogue of $T_{\text{sa}}$ was

$$(5.1) \qquad T_{\text{sv}}(Z) = (Z - \bar{Y})^T \widehat{\Sigma}_Y^{-1}(Z - \bar{Y}) - (Z - \bar{X})^T \widehat{\Sigma}_X^{-1}(Z - \bar{X})$$

and the naive Bayes classifier. Definitions of the two first-mentioned classifiers are given by Chan and Hall (2009). The naive Bayes classifier was constructed under the assumption that all data were normally distributed and employed a ridge parameter. See the last paragraph of of this section for details of the ridging method. In constructing the SV classifier we computed $\widehat{\Sigma}_X$ and $\widehat{\Sigma}_Y$, in (5.1), using the training data from $\Pi_X$ and $\Pi_Y$, respectively, and employing the band-matrix approach studied by Bickel and Levina (2008) with a single band on either side of the main diagonal. Using a single band was appropriate for the small training-sample sizes (3, 5, 8, 15 and 20) encountered with the breast-cancer data.

Training and test datasets were generated and used to assess the five classifiers, as follows. Throughout we took $m = n$. We randomly selected $m$ data vectors from the 623 that represented malignant regions; we similarly chose $n$ from the 101,671 that represented benign regions; we constructed the classifier from these data, and we applied it repeatedly to the remaining $623-m$ data from malignant regions and to a randomly chosen subset of $623-n$ data from the remaining benign data. (Trialling the classifier against all the remaining benign data, i.e., the $101{,}671-n$ benign data not used to construct the classifier, was too time consuming, so we reduced the number to $623-n$, matching that for the malignant data.) This operation was repeated 2000 times, and the error rates averaged to produce the figures discussed below. Note that this procedure gave the two populations prior probabilities of $\frac{1}{2}$ each, rather than the very disparate values of $623/102{,}294 = 0.006$ and $101{,}671/102{,}294 = 0.994$ that would otherwise have prevailed.

Next we summarise the main results. When the common value of $m$ and $n$ was between 15 and 20 the classifiers gave remarkably consistent results over all the settings we treated. In particular, when applied to data from the malignant region the success rates of each of the five classifiers (centroid, SVM, NN, SV and naive Bayes) was in the range 71% to 74%. The ranked order of the classifiers varied from one situation to another, but the centroid-based classifier was almost invariably ranked first. On the other hand (but still for $m$ and $n$ between 15 and 20), when applied to data from the benign region the five classifiers always separated into two clusters on the basis of performance. The centroid and the SV and naive Bayes methods were in the highest-ranked cluster with the centroid method invariably outperforming the naive Bayes approach and the SV method performing close to the centroid method, each having between 72% and 83% success rate. Both of the other two classifiers performed noticeably worse with between 51% and 58% success. Among the latter two methods, either could outperform the other when applied to data from the benign region.

At the extreme of relatively low sample size, and in particular when the common sample size was between 3 and 8, the performance of all classifiers deteriorated and the patterns noted above largely disappeared. For $m$ and $n$ between 5 and 8, and in applications to data from benign regions, the centroid, SV and naive Bayes techniques maintained their superiority over the other two, with the centroid-based method almost invariably the winner. However, in the case of smaller sample sizes the naive Bayes approach had worst performance of all, in both the malignant and benign cases. Here, $m$ and $n$ were far too low for the assumption of normality, on which the naive Bayes method is based, to be even approximately valid. In the case of data from malignant regions the support vector machine also gave good results, being the second best performer behind the centroid method.

Next we give a little more detail in specific cases, starting with the case where $m = 20$. When applied to classify data from malignant regions, the following ranking of classifiers in decreasing order of performance was found: centroid-based method, NN, naive Bayes, SV and SVM. When applied to classify data from benign regions, we found the following rank order: centroid, SV, naive Bayes, SVM and NN. The reasonably good performance of the naive Bayes classifier here was due partly to the fact that when $m = 20$, validity of the assumption of normality was aided by the central limit theorem. In the case of the SV method the larger sample size helped when estimating the covariance matrix. The situation changed markedly when sample sizes were reduced to $m = 5$. There the SV and Bayes methods had significantly more difficulty estimating variance and covariance, to such an extent that using a ridge was essential to obtaining even mediocre performance. When $m = 10$ the Bayes classifier was inferior to each of the other four methods when the data were from malignant regions, and it ranked third, behind the centroid and SV methods, in the case of data from benign regions.

We also explored in more detail the effect of using a ridge parameter to construct the naive Bayes classifier. The ridge was added to conventional estimators of variance, and we sought values of the ridge in the interval $[0.01, 1]$ that maximised classifier success rate, averaged over the malignant and benign cases and for the given choice of $m$. (To put the choice of interval into context we mention that the component-wise average empirical variances of the datasets, for benign and malignant regions, respectively, were 1.00 and 1.21.) Our numerical experiments showed that, when $m = 3$ and the ridge was chosen optimally, the average success rate of the naive Bayes classifier increased from about 50% to 68%. However, when $m = 5$ the average success rate of the naive Bayes classifier increased by only 6%, and the amount of increase declined steadily as $m$ increased; it was only 2% when $m = 20$. Of course, these results are the best possible ones when the ridge is chosen deterministically. In practice the ridge has to be selected empirically, and, especially when $m$ is small (e.g., $m = 3$ or 5), empirical choice of ridge can actually lead to a deterioration in classification performance, since it adds extra noise to the classifier.

### 6. Proof of Theorem 2.

6.1. *Likelihood when* $(C_k, K_k) \equiv (A_k, I_k)$. Let $\phi$ denote the standard normal density. The joint density of $X_{ik}$ (for $1 \leq i \leq m$), $Y_{jk}$ (for $1 \leq j \leq n$) and $Z_k$, for fixed $k$, equals

$$(6.1) \quad E\left[\left\{\prod_{i=1}^{m} \phi(x_{ik} - \delta A_k I_k)\right\}\left\{\prod_{j=1}^{n} \phi(y_{jk} - \delta B_k J_k)\right\} \phi(z_k - \delta C_k K_k)\right]$$
$$= \left\{\prod_{i=1}^{m} \phi(x_{ik})\right\}\left\{\prod_{j=1}^{n} \phi(y_{jk})\right\} \phi(z_k) E(\mathcal{L}_k),$$

where

$$\mathcal{L}_k = \exp\left\{-\frac{1}{2}\delta^2(mA_k^2I_k^2 + nB_k^2J_k^2 + C_k^2K_k^2) \right.$$
$$\left. + \delta\left(A_kI_k\sum_{i=1}^m x_{ik} + B_kJ_k\sum_{j=1}^n y_{jk} + C_kK_kz_k\right)\right\}.$$

Put $\mathcal{X}_k = \{x_{1k},\ldots,x_{mk}\}$, $\mathcal{Y}_k = \{y_{1k},\ldots,y_{nk}\}$, $S_k = \sum_i x_{ik}$ and $T_k = \sum_j y_{jk}$. [Here we keep the data fixed, and so denote them by lower case letters, but from (6.4) down we shall give the data the joint distribution determined by (4.1), and from that point we shall use upper case letters.] If $(C_k, K_k) = (A_k, I_k)$, then

(6.2)
$$E(\mathcal{L}_k \mid \mathcal{X}_k, \mathcal{Y}_k, z_k)$$
$$= E[\exp\{-\tfrac{1}{2}\delta^2(m+1)A_k^2I_k^2 + \delta A_kI_k(S_k+z_k)\} \mid \mathcal{X}_k, \mathcal{Y}_k, z_k]$$
$$\times E\{\exp(-\tfrac{1}{2}\delta^2 nB_k^2J_k^2 + \delta B_kJ_kT_k) \mid \mathcal{X}_k, \mathcal{Y}_k, z_k\}.$$

For $r, s > 0$ and real $t$,

$$E\left\{\exp\left(-\frac{1}{2}r^2sN^2 + rtN\right)\right\} = \exp\left(\frac{1}{2}\frac{r^2t^2}{r^2s+1}\right)(r^2s+1)^{-1/2}.$$

Hence, by (6.2),

(6.3)
$$\psi_1(\mathcal{X}_k, \mathcal{Y}_k, z_k)$$
$$\equiv E(\mathcal{L}_k \mid \mathcal{X}_k, \mathcal{Y}_k, z_k)$$
$$= \left[1 - q + q\{(m+1)\delta^2 + 1\}^{-1/2}\exp\left\{\frac{1}{2}\frac{\delta^2}{(m+1)\delta^2+1}(S_k+z_k)^2\right\}\right]$$
$$\times \left\{1 - q + q(n\delta^2+1)^{-1/2}\exp\left(\frac{1}{2}\frac{\delta^2}{n\delta^2+1}T_k^2\right)\right\}.$$

Combining this result with (6.1) we conclude that the likelihood of $(\mathcal{X}_k, \mathcal{Y}_k, Z_k)$, under the assumption that $(C_k, K_k) = (A_k, I_k)$ is

(6.4)
$$\left\{\prod_{i=1}^m \phi(X_{ik})\right\}\left\{\prod_{j=1}^n \phi(Y_{jk})\right\}\phi(Z_k)\psi_1(\mathcal{X}_k, \mathcal{Y}_k, Z_k).$$

6.2. *Likelihood ratio.* It follows from (6.4) that the ratio of the likelihoods of $(\mathcal{X}_k, \mathcal{Y}_k, Z_k)$, for $(C_k, K_k) = (A_k, I_k)$ versus $(C_k, K_k) = (B_k, J_k)$, is

(6.5)
$$\rho_k(\mathcal{X}_k, \mathcal{Y}_k, Z_k) = \frac{\psi_1(\mathcal{X}_k, \mathcal{Y}_k, Z_k)}{\psi_2(\mathcal{X}_k, \mathcal{Y}_k, Z_k)},$$

where, by symmetry from (6.3),

$$\psi_2(\mathcal{X}_k, \mathcal{Y}_k, Z_k)$$

(6.6)
$$= \left[1 - q + q\{(n+1)\delta^2 + 1\}^{-1/2} \exp\left\{\frac{1}{2} \frac{\delta^2}{(n+1)\delta^2 + 1}(T_k + Z_k)^2\right\}\right]$$

$$\times \left\{1 - q + q(m\delta^2 + 1)^{-1/2} \exp\left(\frac{1}{2} \frac{\delta^2}{m\delta^2 + 1} S_k^2\right)\right\}.$$

The likelihood ratio for the full dataset, $\{(\mathcal{X}_k, \mathcal{Y}_k, Z_k) : 1 \leq k \leq p\}$, is given by

(6.7)
$$\rho = \prod_{k=1}^{p} \rho_k(\mathcal{X}_k, \mathcal{Y}_k, Z_k).$$

6.3. *Properties of $\rho$ when $(C_k, K_k) \equiv (A_k, K_k)$.* Assume that $(C_k, K_k) = (A_k, I_k)$ for all $k$. In this case, writing $N$ for a normal $N(0,1)$ random variable, and interpreting $S_k$, $T_k$ and $Z_k$ as random, we have the following:

$$E\left[\exp\left\{\frac{1}{2} \frac{\delta^2}{(m+1)\delta^2 + 1}(S_k + Z_k)^2\right\}\right]$$

$$= (1-q)E\left[\exp\left\{\frac{1}{2} \frac{(m+1)\delta^2}{(m+1)\delta^2 + 1} N^2\right\}\right]$$

$$+ qE\left(\exp\left[\frac{1}{2} \frac{\{m+1 + (m+1)^2\delta^2\}\delta^2}{(m+1)\delta^2 + 1} N^2\right]\right);$$

$$E\left[\exp\left\{\frac{1}{2} \frac{\delta^2}{(n+1)\delta^2 + 1}(T_k + Z_k)^2\right\}\right]$$

$$= (1-q)^2 E\left[\exp\left\{\frac{1}{2} \frac{(n+1)\delta^2}{(n+1)\delta^2 + 1} N^2\right\}\right]$$

$$+ q^2 E\left(\exp\left[\frac{1}{2} \frac{\{n+1 + (n^2 + 1)\delta^2\}\delta^2}{(n+1)\delta^2 + 1} N^2\right]\right)$$

$$+ q(1-q)\left(E\left[\exp\left\{\frac{1}{2} \frac{(n+1 + n^2\delta^2)\delta^2}{(n+1)\delta^2 + 1} N^2\right\}\right]\right.$$

$$\left. + E\left[\exp\left\{\frac{1}{2} \frac{(n+1 + \delta^2)\delta^2}{(n+1)\delta^2 + 1} N^2\right\}\right]\right);$$

$$E\left\{\exp\left(\frac{1}{2} \frac{\delta^2}{m\delta^2 + 1} S_k^2\right)\right\}$$

$$= (1-q)E\left\{\exp\left(\frac{1}{2} \frac{m\delta^2}{m\delta^2 + 1} N^2\right)\right\} + qE\left[\exp\left\{\frac{1}{2} \frac{(m + m^2\delta^2)\delta^2}{m\delta^2 + 1} N^2\right\}\right];$$

$$E\left\{\exp\left(\frac{1}{2}\frac{\delta^2}{n\delta^2+1}T_k^2\right)\right\}$$
$$=(1-q)E\left\{\exp\left(\frac{1}{2}\frac{n\delta^2}{n\delta^2+1}N^2\right)\right\}+qE\left[\exp\left\{\frac{1}{2}\frac{(n+n^2\delta^2)\delta^2}{n\delta^2+1}N^2\right\}\right].$$

Note too that, for $c<1$, $E\{\exp(\frac{1}{2}cN^2)\}=(1-c)^{-1/2}$. Therefore,

$$E\left[\exp\left\{\frac{1}{2}\frac{\delta^2}{(m+1)\delta^2+1}(S_k+Z_k)^2\right\}\right]$$
$$=\{(m+1)\delta^2+1\}^{1/2}[1-q+q\{1-(m+1)^2\delta^4\}^{-1/2}],$$
$$E\left[\exp\left\{\frac{1}{2}\frac{\delta^2}{(n+1)\delta^2+1}(T_k+Z_k)^2\right\}\right]$$
$$=\{(n+1)\delta^2+1\}^{1/2}((1-q)^2+q^2\{1-(n^2+1)\delta^4\}^{-1/2}$$
$$+q(1-q)\{(1-\delta^4)^{-1/2}+(1-n^2\delta^4)^{-1/2}\}),$$
$$E\left[\exp\left(\frac{1}{2}\frac{\delta^2}{m\delta^2+1}S_k^2\right)\right]=(m\delta^2+1)^{1/2}\{1-q+q(1-m^2\delta^4)^{-1/2}\},$$
$$E\left[\exp\left(\frac{1}{2}\frac{\delta^2}{n\delta^2+1}T_k^2\right)\right]=(n\delta^2+1)^{1/2}\{1-q+q(1-n^2\delta^4)^{-1/2}\}.$$

From these results, (6.3) and (6.6) we see that, if we define

$$\Delta_{S,k}=(m\delta^2+1)^{-1/2}(1-E)\exp\left(\frac{1}{2}\frac{\delta^2}{m\delta^2+1}S_k^2\right),$$
$$\Delta_{T,k}=(n\delta^2+1)^{-1/2}(1-E)\exp\left(\frac{1}{2}\frac{\delta^2}{n\delta^2+1}T_k^2\right),$$
$$\Delta_{SZ,k}=\{(m+1)\delta^2+1\}^{-1/2}(1-E)\exp\left\{\frac{1}{2}\frac{\delta^2}{(m+1)\delta^2+1}(S_k+Z_k)^2\right\},$$
$$\Delta_{TZ,k}=\{(n+1)\delta^2+1\}^{-1/2}(1-E)\exp\left\{\frac{1}{2}\frac{\delta^2}{(n+1)\delta^2+1}(T_k+Z_k)^2\right\},$$
$$\mu_S=\left[(m\delta^2+1)^{-1/2}E\left\{\exp\left(\frac{1}{2}\frac{\delta^2}{m\delta^2+1}S_k^2\right)\right\}-1\right]q^{-1}$$
$$=(1-m^2\delta^4)^{-1/2}-1,$$
$$\mu_T=\left[(n\delta^2+1)^{-1/2}E\left\{\exp\left(\frac{1}{2}\frac{\delta^2}{n\delta^2+1}T_k^2\right)\right\}-1\right]q^{-1}$$
$$=(1-n^2\delta^4)^{-1/2}-1,$$

$$\mu_{SZ} = \left(\{(m+1)\delta^2+1\}^{-1/2} E\left[\exp\left\{\frac{1}{2}\frac{\delta^2}{(m+1)\delta^2+1}(S_k+Z_k)^2\right\}\right] - 1\right) q^{-1}$$

$$= \{1-(m+1)^2\delta^4\}^{-1/2} - 1,$$

$$\mu_{TZ} = \left(\{(n+1)\delta^2+1\}^{-1/2} E\left[\exp\left\{\frac{1}{2}\frac{\delta^2}{(n+1)\delta^2+1}(T_k+Z_k)^2\right\}\right] - 1\right) q^{-1}$$

$$= (1-\delta^4)^{-1/2} + (1-n^2\delta^4)^{-1/2} - 2$$
$$+ q[1+\{1-(n^2+1)\delta^4\}^{-1/2} - (1-\delta^4)^{-1/2} - (1-n^2\delta^4)^{-1/2}],$$

we have,

$$\psi_1(\mathcal{X}_k, \mathcal{Y}_k, Z_k) = (1 + q^2\mu_{SZ} + q\Delta_{SZ,k})(1 + q^2\mu_T + q\Delta_T),$$
$$\psi_2(\mathcal{X}_k, \mathcal{Y}_k, Z_k) = (1 + q^2\mu_{TZ} + q\Delta_{SZ,k})(1 + q^2\mu_S + q\Delta_S).$$

Hence, by (6.5) and (6.7),

$$\rho = \prod_{k=1}^{p} \frac{(1+q^2\mu_{SZ} + q\Delta_{SZ,k})(1+q^2\mu_T + q\Delta_T)}{(1+q^2\mu_{TZ} + q\Delta_{TZ,k})(1+q^2\mu_S + q\Delta_S)}.$$

6.4. *Expansion of likelihood ratio.* Throughout this section we impose the condition, given in Theorem 2, that $\delta = c(mpq^2)^{-1/4}$. The quantities $\mu_S$, $\mu_T$, $\mu_{SZ}$, $\mu_{TZ}$, $\text{var}(\Delta_S)$, $\text{var}(\Delta_T)$, $\text{var}(\Delta_{SZ})$ and $\text{var}(\Delta_{TZ})$, and their counterparts in the case where $(C_k, K_k) \equiv (A_k, K_k)$, are all well defined and finite if and only if, for some $d \in (0, \frac{1}{2})$,

(6.8) $$\max(m+1, n+1)\delta^2 \leq d.$$

This inequality follows from (4.3) and the assumption $\delta = c(mpq^2)^{-1/4}$, provided $c > 0$ is sufficiently small. In this setting we can write

(6.9) $$\rho = \rho_{\text{bias}} \rho_{\text{error}},$$

where

(6.10)
$$\rho_{\text{bias}} = \left\{\frac{(1+q^2\mu_{SZ})(1+q^2\mu_T)}{(1+q^2\mu_{TZ})(1+q^2\mu_S)}\right\}^p,$$

$$\rho_{\text{error}} = \prod_{k=1}^{p} \frac{(1+q\Delta_{SZ,k}/(1+q^2\mu_{SZ}))(1+q\Delta_{T,k}/(1+q^2\mu_T))}{(1+q\Delta_{TZ,k}/(1+q^2\mu_{TZ}))(1+q\Delta_{S,k}/(1+q^2\mu_S))}$$

denote, respectively, the dominant bias term, and the dominant error-about-the-mean term in an expansion of the likelihood ratio $\rho$. We consider two cases:

(i) The ratio $m/n$ is bounded away from zero and infinity as $n \to \infty$.

(a) If $m\delta^2 \to 0$, then

$$1 + \mu_{SZ} = \{1 - (m+1)^2\delta^4\}^{-1/2}$$
$$= (1 - m^2\delta^4)^{-1/2}\left\{1 - \frac{(2m+1)\delta^4}{1 - m^2\delta^4}\right\}^{-1/2}$$
$$= (1 - m^2\delta^4)^{-1/2}\left\{1 + \frac{1}{2}(2m+1)\delta^4 + O(m^3\delta^8)\right\}$$
$$= 1 + \mu_S + \frac{1}{2}(2m+1)\delta^4 + O(m^3\delta^8),$$

(6.11)
$$\mu_{TZ} = \mu_T + \left(1 + \frac{1}{2}\delta^4\right) - 1 + O(\delta^8)$$
$$+ q\left[1 - \left(1 + \frac{1}{2}\delta^4\right)\right.$$
$$\left. + (1 - n^2\delta^4)^{-1/2}\left\{\left(1 - \frac{\delta^4}{1 - n^2\delta^4}\right)^{-1/2} - 1\right\}\right]$$
$$= \mu_T + \frac{1}{2}\delta^4 + O(n^2\delta^8),$$

whence

$$\frac{1 + q^2\mu_{SZ}}{1 + q^2\mu_S} = 1 + \left(m + \frac{1}{2}\right)q^2\delta^4 + O(m^3\delta^8),$$

$$\frac{1 + q^2\mu_T}{1 + q^2\mu_{TZ}} = 1 - \frac{1}{2}q^2\delta^4 + O(q^2n^2\delta^4),$$

(6.12)
$$\rho_{\text{bias}} = \left(\frac{1 + q^2\mu_{SZ}}{1 + q^2\mu_S}\frac{1 + q^2\mu_T}{1 + q^2\mu_{TZ}}\right)^p$$
$$= \exp\{mpq^2\delta^4 + o(mpq^2\delta^4)\}.$$

To treat $\rho_{\text{error}}$, note that

$$\{(m+1)\delta^2 + 1\}^{-1/2} = (m\delta^2 + 1)^{-1/2}\{1 + O(\delta^4)\},$$

$$\exp\left\{\frac{1}{2}\frac{\delta^2}{(m+1)\delta^2 + 1}(S_k + Z_k)^2\right\} = \exp\left\{\frac{1}{2}\frac{\delta^2}{m\delta^2 + 1}(S_k + Z_k)^2\right\}$$
$$\times (1 + m\delta^4 R_k),$$

where, here and below, $R_1, R_2, \ldots$ is a generic sequence of independent and identically distributed random variables, depending on $\delta$ but for which, for each $r \geq 1$, absolute moments of order $r$ are uniformly bounded provided $\delta$ is sufficiently small. Therefore, recalling the genericity of the notation $R_k$,

$$\{(m+1)\delta^2 + 1\}^{-1/2} \exp\left\{\frac{1}{2}\frac{\delta^2}{(m+1)\delta^2+1}(S_k+Z_k)^2\right\}$$

$$= (m\delta^2+1)^{-1/2}\exp\left\{\frac{1}{2}\frac{\delta^2}{m\delta^2+1}(S_k+Z_k)^2\right\}(1+m\delta^4 R_k)$$

$$= (m\delta^2+1)^{-1/2}\exp\left(\frac{1}{2}\frac{\delta^2}{m\delta^2+1}S_k^2\right)$$

$$\times \left\{1 + \frac{1}{2}\frac{\delta^2}{m\delta^2+1}(2S_k Z_k + Z_k^2) + m\delta^4 R_k\right\}.$$

Hence,

$$\Delta_{SZ,k} = \Delta_{S,k} + (1-E)(m\delta^2+1)^{-1/2}\exp\left(\frac{1}{2}\frac{\delta^2}{m\delta^2+1}S_k^2\right)$$

$$\times \left\{\frac{1}{2}\frac{\delta^2}{m\delta^2+1}(2S_k Z_k + Z_k^2)\right\} + (1-E)m\delta^4 R_k,$$

and from (6.11), $\mu_{SZ} = \mu_S + O(m\delta^4)$. Therefore,

$$1 + q\frac{\Delta_{SZ,k}}{1+q^2\mu_{SZ}} = 1 + q\frac{\Delta_{S,k}}{1+q^2\mu_S}$$

$$+ q(1-E)(m\delta^2+1)^{-1/2}\exp\left(\frac{1}{2}\frac{\delta^2}{m\delta^2+1}S_k^2\right)$$

$$\times \left\{\frac{1}{2}\frac{\delta^2}{m\delta^2+1}(2S_k Z_k + Z_k^2)\right\} + (1-E)mq\delta^4 R_k,$$

whence, since $\Delta_{S,k} = (1-E)m\delta^2 R_k$,

(6.13) $\quad \dfrac{1+q\Delta_{SZ,k}/(1+q^2\mu_{SZ})}{1+q\Delta_{S,k}/(1+q^2\mu_S)} = 1 + U_k + (1-E)mq\delta^4 R_k,$

where

$$U_k = q(1-E)(m\delta^2+1)^{-1/2}\exp\left(\frac{1}{2}\frac{\delta^2}{m\delta^2+1}S_k^2\right)$$

$$\times \left\{\frac{1}{2}\frac{\delta^2}{m\delta^2+1}(2S_k Z_k + Z_k^2)\right\}.$$

Analogously,

(6.14) $\quad \dfrac{1+q\Delta_{TZ,k}/(1+q^2\mu_{TZ})}{1+q\Delta_{T,k}/(1+q^2\mu_T)} = 1 + V_k + (1-E)nq\delta^4 R_k,$

where

$$V_k = 1 + q(1-E)(n\delta^2+1)^{-1/2}\exp\left(\frac{1}{2}\frac{\delta^2}{n\delta^2+1}T_k^2\right)$$
$$\times\left\{\frac{1}{2}\frac{\delta^2}{n\delta^2+1}(2T_kZ_k + Z_k^2)\right\}.$$

Since the assumption that $m\delta^2$ is bounded entails $p^{1/2}mq\delta^4 = o(m^{1/2}p^{1/2}q\delta^2)$, then (6.10), (6.13) and (6.14) imply that

$$\rho_{\text{error}} = \prod_{k=1}^{p}\frac{(1+q\Delta_{SZ,k}/(1+q^2\mu_{SZ}))(1+q\Delta_{T,k}/(1+q^2\mu_T))}{(1+q\Delta_{TZ,k}/(1+q^2\mu_{TZ}))(1+q\Delta_{S,k}/(1+q^2\mu_S))}$$

(6.15)
$$= \{1 + O_p(p^{1/2}mq\delta^4)\}\prod_{k=1}^{p}\{(1+U_k)/(1+V_k)\}$$
$$= \exp\left\{\sum_{k=1}^{p}(U_k - V_k) - \frac{1}{2}\sum_{k=1}^{p}(U_k^2 - V_k^2) + o_p(1)\right\}.$$

Now, $W = \sum_k(U_k - V_k)$ is asymptotically normal $N\{0, (m+n)pq^2\delta^4\}$, and $\sum_k(U_k^2 - V_k^2) = (m-n)pq^2\delta^4 + o_p(1)$. These properties and (6.15) imply that

(6.16) $$\rho_{\text{error}} = \exp\{W - \tfrac{1}{2}(m-n)pq^2\delta^4 + o_p(1)\}.$$

Combining (6.9), (6.12) and (6.16) we deduce that

(6.17) $$\rho = \exp[N\{(m+n)pq^2\delta^4\}^{1/2} + \tfrac{1}{2}(m+n)pq^2\delta^4 + o_p(1)],$$

where $N$ is asymptotically normal $N(0,1)$. Therefore, if $\widehat{\chi}$ is taken to be the likelihood-ratio classifier then for all values of $c$ that are sufficiently small to ensure that (6.8) holds for some $d < \tfrac{1}{2}$, then

(6.18) $$\liminf_{n\to\infty}\{P_X(\widehat{\chi}=B) + P_Y(\widehat{\chi}=A)\} > 0.$$

This establishes Theorem 2 in the case where $m\delta^2 \to 0$.

(b) If $\ell_1 \equiv m\delta^2$ and $\ell_2 \equiv n\delta^2 \to$ converge to finite, nonzero constants, both of them strictly less than 1, then

$$1 + \mu_{SZ} = \{1 - (m+1)^2\delta^4\}^{-1/2}$$
$$= (1-\ell_1^2)^{-1/2}\left\{1 - \frac{(2m+1)\delta^4}{1-\ell_1^2}\right\}^{-1/2}$$
$$= 1 + \mu_S + (1-\ell_1^2)^{-3/2}\left(m+\frac{1}{2}\right)\delta^4 + O(m^2\delta^8),$$

$$\mu_{TZ} = \mu_T + \left(1 + \frac{1}{2}\delta^4\right) - 1 + O(\delta^8)$$
$$+ q\left[1 - \left(1 + \frac{1}{2}\delta^4\right) + (1 - n^2\delta^4)^{-1/2}\left\{\left(1 - \frac{\delta^4}{1 - n^2\delta^4}\right)^{-1/2} - 1\right\}\right]$$
$$= \mu_T + \frac{1}{2}\delta^4[1 + q\{(1 - \ell_2^2)^{-3/2} - 1\}] + O(\delta^8),$$

whence

$$\frac{1 + q^2\mu_{SZ}}{1 + q^2\mu_S} = 1 + q^2\delta^4\frac{(1 - \ell_1^2)^{-3/2}(m + 1/2)}{1 + q^2\{(1 - \ell_1^2)^{-1/2} - 1\}}$$
$$+ O(q^2 m^2 \delta^8),$$
$$\frac{1 + q^2\mu_T}{1 + q^2\mu_{TZ}} = 1 - \frac{1}{2}q^2\delta^4\frac{1 + q\{(1 - \ell_2^2)^{-3/2} - 1\}}{1 + q^2\{(1 - \ell_2^2)^{-1/2} - 1\}} + O(q^2\delta^8),$$
$$\rho_{\text{bias}} = \left\{\frac{1 + q^2\mu_{SZ}}{1 + q^2\mu_{TZ}}\frac{1 + q^2\mu_T}{1 + q^2\mu_S}\right\}^p$$
$$= \exp\{L_1 m p q^2 \delta^4 + o(1)\},$$

where

$$L_j = \frac{(1 - \ell_j^2)^{-3/2}}{1 + q^2\{(1 - \ell_j^2)^{-1/2} - 1\}}.$$

Compare (6.12). A similar argument can be used to derive an analogue of (6.15) in this setting, giving, via (6.9), the following analogue of (6.17):

$$\rho = \exp[N\{(L_1 m + L_2 n)pq^2\delta^4\}^{1/2} + \tfrac{1}{2}(L_1 m + L_2 n)pq^2\delta^4 + o_p(1)],$$

where $N$ is asymptotically normal $N(0,1)$. Result (6.18) follows as before.

DEPARTMENT OF MATHEMATICS AND STATISTICS
UNIVERSITY OF MELBOURNE
PARKVILLE, VIC, 3010
AUSTRALIA
E-MAILS: halpstat@ms.unimelb.edu.au
tungp@ms.unimelb.edu.au